\documentclass[12pt,twoside]{article}
\usepackage{amsmath}
\usepackage{amssymb}
\usepackage{graphicx}
\def\eqvsp{}  \newdimen\paravsp  \paravsp=1.3ex
\def\,{\mskip 3mu} \def\>{\mskip 4mu plus 2mu minus 4mu} \def\;{\mskip 5mu plus 5mu} \def\!{\mskip-3mu}
\def\dispmuskip{\thinmuskip= 3mu plus 0mu minus 2mu \medmuskip=  4mu plus 2mu minus 2mu \thickmuskip=5mu plus 5mu minus 2mu}
\def\textmuskip{\thinmuskip= 0mu                    \medmuskip=  1mu plus 1mu minus 1mu \thickmuskip=2mu plus 3mu minus 1mu}
\textmuskip
\def\beq{\eqvsp\dispmuskip\begin{equation}}    \def\eeq{\eqvsp\end{equation}\textmuskip}
\def\beqn{\eqvsp\dispmuskip\begin{displaymath}}\def\eeqn{\eqvsp\end{displaymath}\textmuskip}
\def\bqa{\eqvsp\dispmuskip\begin{eqnarray}}    \def\eqa{\eqvsp\end{eqnarray}\textmuskip}
\def\bqan{\eqvsp\dispmuskip\begin{eqnarray*}}  \def\eqan{\eqvsp\end{eqnarray*}\textmuskip}

\newenvironment{keyword}{\centerline{\bf\small
Keywords}\begin{quote}\small}{\par\end{quote}\vskip 1ex}
\newtheorem{theorem}{Theorem}

\newtheorem{myexample}[theorem]{Example}
\def\fexample#1#2#3{\vspace{-1ex}\begin{myexample}[#2]\label{#1}\rm #3
\hspace*{\fill} $\diamondsuit\quad$ \end{myexample}\vspace{-1ex} }
\def\paradot#1{\vspace{\paravsp plus 0.5\paravsp minus 0.5\paravsp}\noindent{\bf\boldmath{#1.}}}

\def\req#1{(\ref{#1})}

\def\eps{\varepsilon}
\def\nq{\hspace{-1em}}
\def\qed{\hspace*{\fill}\rule{1.4ex}{1.4ex}$\quad$\\}
\def\fr#1#2{{\textstyle{#1\over#2}}}
\def\frs#1#2{{^{#1}\!/\!_{#2}}}
\def\SetR{I\!\!R}
\def\SetN{I\!\!N}

\def\qmbox#1{{\quad\mbox{#1}\quad}}

\def\v{\boldsymbol}

\def\s{\sigma}

\def\p{{\scriptscriptstyle+}}
\def\pp{{\scriptscriptstyle++}}
\def\n{{n}}

\def\u{u} 

\def\v{\boldsymbol} 
\def\vt{{\v t}}\def\vu{{\v u}}\def\vpi{{\v\pi}}
\def\pin{{\scriptstyle\Pi}}
\def\Var{{\mbox{Var}}}

\def\propersubset{\subsetneq}
\def\low{\underline}
\def\up{\overline}
\def\DeltaOX{\Delta\hspace{-8pt}{\scriptscriptstyle^{_\otimes}}}
\def\deltaOX{\Delta\hspace{-6.5pt}{\scriptscriptstyle^{_\otimes}}}
\def\Deltapi{\Delta}
\def\Deltal{\Delta_e}
\def\Deltapl{\Delta'_e}
\def\ots{{\scriptscriptstyle\!\otimes}}
\def\Vol{\mbox{Vol}}
\def\leqsq{\sqsubseteq}
\def\geqsq{\sqsupseteq}
\def\maxo{\max}                        
\def\mino{\min}                        
\def\indfct{{1\!\!1\!}}

\begin{document}

\title{\vspace{-4ex}
\vskip 2mm\bf\Large\hrule height5pt \vskip 4mm
Practical Robust Estimators for \\ the Imprecise Dirichlet Model
\vskip 4mm \hrule height2pt}
\author{{\bf Marcus Hutter}\\[3mm]
\normalsize RSISE$\,$@$\,$ANU and SML$\,$@$\,$NICTA \\
\normalsize Canberra, ACT, 0200, Australia \\
\normalsize \texttt{marcus@hutter1.net \ \  www.hutter1.net}%
\thanks{A shorter version appeared in the proceedings of the ISIPTA 2003 conference \cite{Hutter:03idm}.}
}
\date{January 2009}
\maketitle

\begin{abstract}
Walley's Imprecise Dirichlet Model (IDM) for categorical i.i.d.\
data extends the classical Dirichlet model to a set of priors. It
overcomes several fundamental problems which other approaches to
uncertainty suffer from. Yet, to be useful in practice, one needs
efficient ways for computing the imprecise=robust sets or
intervals. The main objective of this work is to derive exact,
conservative, and approximate, robust and credible interval
estimates under the IDM for a large class of statistical
estimators, including the entropy and mutual information.
\vspace{3ex}
\def\contentsname{\centering\normalsize Contents}
{\parskip=-2.5ex\tableofcontents}
\end{abstract}


\hfill\newpage
\section{Introduction}\label{secInt}

This work derives interval estimates under the Imprecise Dirichlet
Model (IDM) \cite{Walley:96} for a large class of statistical
estimators. In the IDM one considers an i.i.d.\ process with unknown
chances\footnote{Also called {\em objective} or {\em aleatory}
probabilities.} $\pi_i$ for outcome $i\in\{1,...,d\}$. The prior
uncertainty about\footnote{We denote vectors by $\v
x:=(x_1,...,x_d)$ for $\v x\in\{\v n,\vt,\vu,\vpi,...\}$, and $i$
ranges from $1$ to $d$ unless otherwise stated. See also Appendix
\ref{secSymb}.} $\vpi=(\pi_1,...,\pi_d)$ is modeled by a set of
Dirichlet priors\footnote{Also called {\em second
order} or {\em subjective} or {\em belief} or {\em epistemic}
probabilities.} $\{p(\vpi)\propto\prod_i\pi_i^{st_i-1}\,:\,\vt\in\Delta\}$, where%
\footnote{Strictly speaking, $\Delta$ should be the open simplex
\cite{Walley:96}, since $p(\vpi)$ is improper for $\vt$ on the
boundary of $\Delta$. For simplicity we assume that, if necessary,
considered functions of $\vt$ can and are continuously extended to
the boundary of $\Delta$, so that, for instance, minima and maxima
exist. All considerations can straightforwardly, but cumbersomely,
be rewritten in terms of an open simplex. Note that open/closed
$\Delta$ result in open/closed robust intervals, the difference
being numerically/practically irrelevant.} %
$\Delta:=\{\vt\,:\,t_i\geq 0\,\forall i,\, \sum_i t_i=1\}$, and $s$
is a hyper-parameter, typically chosen between 1 and 2. Sets of
probability distributions are often called Imprecise probabilities,
hence the name IDM for this model. We avoid the term {\em imprecise}
and use {\em robust} instead, or capitalize {\em Imprecise}. The IDM
overcomes several fundamental problems which other approaches to
uncertainty suffer from \cite{Walley:96}. For instance, the IDM
satisfies the representation invariance principle and the symmetry
principle, which are mutually exclusive in a pure Bayesian treatment
with proper prior \cite{Walley:96}.$\!\!$\footnote{ But see
\cite{Hutter:07uspx} for a proper Bayesian reconciliation of these
principles.} The counts $n_i$ for $i$ form a minimal sufficient
statistic of the data of size $n=\sum_i n_i$. Statistical estimators
$F(\v n)$ usually also depend on the chosen prior: so a set of
priors leads to a set of estimators $\{F_\vt(\v
n)\,:\;\vt\in\Delta\}$. For instance, the expected chances
$E_\vt[\pi_i]={n_i+st_i\over n+s}=:\u_i(\vt)$ lead to a robust
interval estimate $[{n_i\over n+s},{n_i+s\over n+s}]\ni
E_\vt[\pi_i]$. Robust intervals for the variance $\Var_\vt[\pi_i]$
\cite{Walley:96} and for the mean and variance of
linear-combinations $\sum_i\alpha_i\pi_i$ have also been derived
\cite{Bernard:01}. Bayesian estimators (like expectations) depend on
$\vt$ and $\v n$ only through $\vu$ (and $n+s$ which we suppress),
i.e.\ $F_\vt(\v n)=F(\vu)$. The main objective of this work is to
derive approximate, conservative, and exact intervals
$[\mino_{\vt\in\Delta}F(\vu),\maxo_{\vt\in\Delta}F(\vu)]$ for
general $F(\vu)$, and for the expected (also called predictive)
entropy and the expected mutual information in particular. These
results are key building blocks for applying the IDM. Walley
suggests, for instance, to use $\mino_\vt P_\vt[{\cal F}\geq
c]\geq\alpha$ for inference problems and $\mino_\vt E_\vt[{\cal
F}]\geq c$ for decision problems \cite{Walley:96}, where $\cal F$ is
some function of $\vpi$. One application is the inference of robust
tree-dependency structures \cite{Zaffalon:01tree,Hutter:05tree}, in
which edges are partially ordered based on Imprecise mutual
information.

Section \ref{secIDM} gives a brief introduction to
the IDM and describes our problem setup.
In Section \ref{secEEI} we derive exact robust intervals for
concave functions $F$, such as the entropy.
Section \ref{secAEI} derives approximate robust intervals for arbitrary
$F$.
In Section \ref{secEP} we show how bounds of elementary functions
can be used to get bounds for composite function, especially for
sums and products of functions. The results are used in
Section \ref{secIEMI} for deriving robust intervals for the
mutual information.
The issue of how to set up IDM models on
product spaces is discussed in Section \ref{secPS}.
Section \ref{secCI} addresses the problem of how to combine
Bayesian credible intervals with the robust intervals of the IDM.
Conclusions are given in Section \ref{secConc}.
Appendix \ref{secPsi} lists properties of the $\psi$ function,
which occurs in the expressions for the expected entropy and
mutual information.
Appendix \ref{secSymb} contains a table of used notation.

\section{The Imprecise Dirichlet Model}\label{secIDM}

This section provides a brief introduction to the IDM, introduces
notation, and describes our generic problem setup of finding upper
and lower statistical estimators. We first introduce the multinomial
process and the Bayesian treatment with Dirichlet priors, and then
the IDM extension to sets of such priors. See \cite{Walley:96} for a
more thorough account and motivation.

\paradot{Random i.i.d.\ processes}
We consider discrete random variables $\imath\in\{1,...,d\}$ and an
i.i.d.\ random process with outcome $i\in\{1,...,d\}$ having
probability $\pi_i$. The chances $\vpi$ form a probability
distribution, i.e.\ $\vpi\in\Deltapi:=\{\v x\in\SetR^d\,:\,x_i\geq
0\,\forall i,\; x_\p=1\}$, where we have used the abbreviation $\v
x=(x_1,...,x_d)$ and $x_\p:=\sum_{i=1}^d x_i$. The likelihood of a
specific (ordered) data set $\v D=(i_1,...,i_n)$ with $n_i$
observations $i$ and total sample size $n=n_\p=\sum_i n_i$ is $p(\v
D|\vpi)=\prod_i\pi_i^{n_i}$.
The chances $\pi_i$ are usually unknown and have to be estimated
from the sample frequencies $n_i$. The maximum likelihood
(frequency) estimate ${n_i\over n}$ for $\pi_i$ is one possible
point estimate.

\paradot{The Bayesian approach}
A (precise) Bayesian models the initial uncertainty in $\vpi$ by a
(second order) prior ``belief'' distribution $p(\vpi)$ with domain
$\vpi\in\Deltapi$. The Dirichlet priors
$p(\vpi)\propto\prod_i\pi_i^{n'_i-1}$, where $n'_i$ comprises prior
information, represent a large class of priors. The $n'_i$ may be
interpreted as (possibly fractional) virtual number of
``observations''. High prior belief in $i$ can be modeled by large
$n'_i$. It is convenient to write $n'_i=s\cdot t_i$ with $s:=n'_+$,
hence $\vt\in\Delta$. Having no initial bias one should choose a
prior in which all $t_i$ are equal, i.e.\ $t_i={1\over d}\,\forall
i$. Examples for $s$ are $0$ for Haldane's prior \cite{Haldane:48},
$1$ for Perks' prior \cite{Perks:47}, ${d\over 2}$ for Jeffreys'
prior \cite{Jeffreys:46}, and $d$ for Bayes-Laplace's uniform prior
\cite{Gelman:95}. From the prior and the data likelihood one can
determine the posterior $p(\vpi|\v D)=p(\vpi|\v
n)\propto\prod_i\pi_i^{n_i+st_i-1}$.

The posterior $p(\vpi|\v D)$ summarizes all statistical
information available in the data. In general, the posterior is a
very complex object, so we are interested in summaries of this
plethora of information. A possible summary is the expected
value or mean
$
  E_\vt[\pi_i]={n_i+st_i\over n+s}
$
which is often used for estimating $\pi_i$. The accuracy may be
obtained from the covariance of $\vpi$.

Usually one is not only interested in an estimation of the whole
vector $\vpi$, but also in an estimation of scalar functions ${\cal
F}:\Deltapi\to\SetR$ of $\vpi$, such as the entropy ${\cal
H}(\vpi)=-\sum_i\pi_i\log\pi_i$, where $\log$ denotes the natural
logarithm. Since $\cal F$ is itself a random variable we could
determine the posterior distribution $p({\cal F}_0|\v
n)=\int_{\Deltapi}\delta({\cal F}(\vpi)-{\cal F}_0)p(\vpi|\v
n)d\vpi$ of ${\cal F}$, where ${\cal F}_0\in\SetR$ and $\delta()$ is
the Dirac delta distribution. This may further be summarized by the
posterior mean $E_\vt[{\cal F}]=\int_{\Deltapi}{\cal
F}(\vpi)p(\vpi|\v n)d\vpi$ and possibly the posterior variance
$\Var_\vt[{\cal F}]$.
A simple but crude approximation for the mean can be obtained by
exchanging $E$ with ${\cal F}$ (exact only for linear functions):
$E_\vt[{\cal F}(\vpi)]\approx {\cal F}(E_\vt[\vpi])$. The
approximation error is typically of the order $1\over n$.

\paradot{The Imprecise Dirichlet Model}
There are several problems with this approach. First, the uniform
choice $t_i={1\over d}$ depends on {\em how} events are grouped
into $d$ classes, which could be ambiguous.
Secondly, it assumes exact prior knowledge of $p(\vpi)$. The
solution to the second problem is to model our ignorance by
considering sets of priors $p(\vpi)$, often called Imprecise
probabilities. The specific {\em Imprecise Dirichlet Model} (IDM)
\cite{Walley:96} considers the set of {\em all} $\vt\in\Delta$,
i.e. $\{p(\vpi|\v n):\vt\in\Delta\}$ which solves also the first
problem. Walley suggests to fix the hyperparameter $s$ somewhere
in the interval $[1,2]$. A set of priors results in a set of
posteriors, set of expected values, etc. For real-valued
quantities like the expected entropy $E_\vt[{\cal H}]$ the sets
are typically intervals, which we call robust
intervals
\beqn
  E_\vt[{\cal F}] \in [\mino_{\vt\in\Delta}E_\vt[{\cal F}] \,,\,
    \maxo_{\vt\in\Delta}E_\vt[{\cal F}]].
\eeqn

\paradot{Problem setup and notation}
Consider any statistical estimator $F$. $F$ is a function of the
data $\v D$ and the hyperparameters $\vt$.
We define the general correspondence
\beq\label{thtrel}
 \u_i^{\cdots}={n_i+st_i^{\cdots}\over n+s},
 \quad\mbox{where $^{\ldots}$ can be various superscripts or be empty}.
\eeq
$F$ can, hence, be rewritten as a function of $\vu$ and $\v D$.
Since we regard $\v D$ as fixed, we suppress this dependence and
simply write $F=F(\vu)$. This is further motivated by the fact that
all Bayesian estimators of functions $\cal F$ of $\vpi$ only depend
on $\vu$ and the sample size $n+s$. It is easy to see that this
holds for the mean, i.e.\ $E_\vt[{\cal F}]=F(\vu\,;\,n+s)$, and
similarly for the variance and all higher (central) moments. Most of
this work is applicable to generic $F$, whatever it's origin -- as
an expectation of $\cal F$ or otherwise. The main focus of this work
is to derive exact and approximate expressions for upper and lower
$F$ values
\beqn
 \up{F}:=\maxo_{\vt\in\Delta}F(\vu)
 \qmbox{and}
 \low F:=\mino_{\vt\in\Delta}F(\vu),\qquad
 \up{\low F}:=[\low F,\up F].
\eeqn
$\vt\in\Delta$ $\Leftrightarrow$
$\vu\in\Delta'$, where $\Delta':=\{\vu\,:\,\u_i\geq{n_i\over
n+s}\,\forall i,\;\u_\p=1\}$. We define $\vu^{\up F}$ as the
$\vu\in\Delta'$ which maximizes $F$, i.e.\ $\up F=F(\vu^{\up
F})$, and similarly $\vt^{\up F}$ through relation
\req{thtrel}. If the maximum of $F$ is assumed in a corner of
$\Delta'$ we denote the index of the corner by $i^{\up F}$, i.e.\
$t_i^{\up F}=\delta_{ii^{\up F}}$, where $\delta_{ij}$ is
Kronecker's delta function, and similarly for
$\vu^{\low F}$, $\vt^{\low F}$, $i^{\low F}$.

\section{Exact Robust Intervals for Concave Estimators}\label{secEEI}

In this section we derive exact expressions for $\up{\low F}$ if
$F:\Delta\to\SetR$ is of the form
\beq\label{Fconc}
  F(\vu)=\sum_{i=1}^d f(\u_i)
  \qmbox{and concave} f:[0,1]\to\SetR.
\eeq
The expected entropy is such an example (discussed later). Convex
$f$ are treated similarly (or simply take $-f$).

\paradot{The nature of the solution}
The approach to a solution of this problem is motivated as follows:
Due to symmetry and concavity of $F$, the global maximum is attained
at the center $\u_i={1\over d}$ of the probability simplex $\Delta$
if we allow $\vu\in\Delta$, i.e.\ the more uniform $\vu$ is, the
larger $F(\vu)$. The nearer $\vu$ is to a vertex of $\Delta$, i.e.\
the more unbalanced $\vu$ is, the smaller is $F(\vu)$. But the
constraints $t_i\geq 0$ restrict $\vu$ to the smaller simplex
\beqn
  \Delta'=\{\vu\,:\,\u_i\geq\u_i^0\,\forall i,\;\u_\p=1\}
  \qmbox{with} \u_i^0:={n_i\over n+s},
\eeqn
which prevents setting $\u_i^{\up
F}={1\over d}$ and $\u_i^{\low F}=\delta_{i1}$.
Nevertheless, the basic idea of choosing $\vu$ as
uniform / as unbalanced as possible still works, as we will see.

\paradot{Greedy $F(\vu)$ minimization}
Consider the following procedure for obtaining $\vu^{\low F}$. We
start with $\vt\equiv\v 0$ (outside the usual domain $\Delta$ of
$F$, which can be extended to $[0,1]^d$ via \req{Fconc}) and then
gradually increase $\vt$ in an axis-parallel way until $t_\p=1$.
With axis-parallel we mean that only one component of $\vt$ is
increased, which one possibly changes during the process. The total
zigzag curve from $\vt^{start}=\v 0$ to $\vt^{end}$ has length
$t_\p^{end}=1$. Since all possible curves have the same (Manhattan)
length 1, $F(\vu^{end})$ is minimized for the curve which has (on
average) smallest $F$-gradient along its path. A greedy strategy is
to follow the direction $i$ of currently smallest $F$-gradient
${\partial F\over\partial t_i}=f'(\u_i){s\over n+s}$. Since $f'$ is
monotone decreasing ($f''<0$), ${\partial F\over\partial t_i}$ is
smallest for largest $\u_i$. At $\vt^{start}=\v 0$, $\u_i={n_i\over
n+s}$ is largest for $i=i^{min}:=\arg\max_i n_i$. Once we start in
direction $i^{min}$, $\u_{i^{min}}$ increases even further whereas
all other $\u_i$ ($i\neq i^{min}$) remain constant. So the moving
direction is never changed and finally we reach a local minimum at
$t_i^{end}=\delta_{ii^{min}}$. Below we show that this is a global
minimum, i.e.
\beq\label{Fmin}
  t_i^{\low F}=\delta_{ii^{\low F}} \qmbox{with}
  i^{\low F}:=\arg\max_i n_i.
\eeq

\paradot{Greedy $F(\vu)$ maximization}
Similarly we maximize $F(\vu)$. Now we increase $\vt$ in direction
$i=i_1$ of maximal ${\partial F\over\partial t_i}$, which is
the direction of smallest $\u_i$. 
Again, (only) $\u_{i_1}$ increases, but possibly reaches a value
where it is no longer the smallest one. We stop if it becomes equal
to the second smallest $\u_i$, say $i=i_2$. We now have to increase
$\u_{i_1}$ and $\u_{i_2}$ with same speed (or in an $\eps$-zigzag
fashion) until they become equal to $\u_{i_3}$, etc.\ or until
$\u_\p=1=t_\p$ is reached. Assume the process stops with direction
$i_m$ and minimal $\u$ being $\tilde\u$, i.e.\ finally
$\u_{i_k}=\tilde\u$ for $k\leq m$ and $t_{i_k}=0$ for $k>m$. From
the constraint $1=\u_\p=\sum_{k\leq m}\u_{i_k}+\sum_{k>m}\u_{i_k} =
m\tilde\u+\sum_{k>m}{n_{i_k}\over n+s}$ we obtain $\tilde\u={1\over
m}[1-\sum_{k>m}{n_{i_k}\over n+s}]= [s+\sum_{k\leq
m}n_{i_k}]/[m(n+s)]$. One can show that $\tilde\u$ as a function of
$m$ has one global minimum (no local ones) and that the final $m$ is
the one which minimizes $\tilde\u$, i.e.\
\beq\label{Fmax}
  \tilde\u = \min_{m\in\{1...d\}} {s+\sum_{k\leq m} n_{i_k}\over
  m(n+s)}, \;\;\mbox{where}\;\; n_{i_1}\leq
  n_{i_2}\leq...\leq\n_{i_d},\quad
  \u_i^{\up F}=\max\{\u_i^0,\tilde\u\}.
\eeq
If there is a unique minimal $n_{i_1}$ with gap $\geq s$ to the 2nd
smallest $n_{i_2}$ (which is quite likely for not too small $n$ and
small $s$ like 1 or 2), then $m=1$ and the maximum is attained at a
corner of $\Delta$ ($\Delta'$).

\begin{theorem}[Exact extrema for concave functions on simplices]\label{thEEI}
Assume $F:\Delta'\to\SetR$ is a concave function of the form
$F(\vu)=\sum_{i=1}^df(\u_i)$.
Then $F$ attains the global maximum $\up F$ at $\vu^{\up F}$
defined in \req{Fmax} and the global minimum $\low F$ at
$\vu^{\low F}$ defined in \req{Fmin}.
\end{theorem}

\paradot{Proof} What remains to be shown is that the solutions
obtained in the last paragraphs by greedy
minimization/maximization of $F(\vu)$ are actually global
minima/maxima. For this assume that $\vt$ is a local minimum of
$F(\vu)$. Let $j:=\arg\max_i\u_i$ (ties broken arbitrarily).
Assume that there is a $k\neq j$ with non-zero $t_k$. Define
$\vt'$ as $t'_i=t_i$ for all $i\neq j,k$, and $t'_j=t_j+\eps$,
$t'_k=t_k-\eps$, for some $0<\eps\leq t_k$. From $\u_k\leq \u_j$
and the concavity of $f$ we get\footnote{Slope
${f(\u+\eps)-f(\u)\over\eps}$ is a decreasing function in $\u$
for any $\eps>0$, since $f$ is concave.}
\bqan
  F(\vu')-F(\vu) &=& [f(\u'_j)+f(\u'_k)]-[f(\u_j)+f(\u_k)]
  \\ &=&
  [f(\u_j\!+\!\sigma\eps)-f(\u_j)]-[f(\u_k)-f(\u_k\!-\!\sigma\eps)] \;<\; 0,
\eqan
where $\sigma:={s\over n+s}$. This contradicts the minimality
assumption of $\vt$. Hence, $t_i=0$ for all $i$ except one (namely $j$,
where it must be 1). (Local) minima are attained in the vertices
of $\Delta$. Obviously the global minimum is for $t_i^{\low
F}=\delta_{ii^{\low F}}$ with $i^{\low F}:=\arg\max_i n_i$. This
solution coincides with the greedy solution. Note that the global
minimum may not be unique, but since we are only interested in
the value of $F(\vu^{\low F})$ and not its argument this
degeneracy is of no further significance.

Similarly for the maximum, assume that $\vt$ is a (local)
maximum of $F(\vu)$. Let $j:=\arg\min_i\u_i$ (ties broken
arbitrarily). Assume that there is a $k\neq j$ with
non-zero $t_k$ {\em and} $\u_k>\u_j$. Define $\vt'$ as above with
$0<\eps<\min\{t_k\,,\,t_k-t_j\}$. Concavity of $f$ implies
\beqn
  F(\vu')-F(\vu) =
  [f(\u_j\!+\!\sigma\eps)-f(\u_j)]-[f(\u_k)-f(\u_k\!-\!\sigma\eps)] > 0,
\eeqn
which contradicts the maximality assumption of $\vt$. Hence
$t_i=0$ if $\u_i$ is not minimal ($\tilde\u$). The previous
paragraph constructed the unique solution $\vu^{\up F}$ satisfying
this condition. Since this is the only local maximum it must be
the unique global maximum (contrast this to the minimum case).
\qed

\begin{theorem}[Exact extrema of expected entropy]\label{corEEI}
Let ${\cal H}(\vpi)=-\sum_i \pi_i\log\pi_i$ be the entropy of
$\vpi$ and the uncertainty of $\vpi$ be modeled by the Imprecise
Dirichlet Model. The expected entropy $H(\vu):=E_\vt[{\cal
H}]$ for given hyperparameter $\vt$ and sample $\v n$ is given by
\beq\label{hex}
  H(\vu)=\sum_i h(\u_i) \qmbox{with}
  h(\u)=\u\!\cdot\![\psi(n\!+\!s\!+\!1)-\psi((n\!+\!s)\u\!+\!1)] =
  \u\cdot\;\nq\nq\sum_{k={(n+s)}u+1}^{n+s}\nq\!\! k^{-1}
\eeq
where $\psi(x)=d\,\log\Gamma(x)/dx$ is the logarithmic derivative of
the Gamma function and the last expression is valid for integral $s$
and $(n+s)u$. The lower $\low H$ and upper $\up H$ expected
entropies are assumed at $\vu^{\low H}$ and $\vu^{\up H}$ given in
\req{Fmin} and \req{Fmax} (with $F$ replaced by $H$, see also
\req{thtrel}).
\end{theorem}

A derivation of the exact expression \req{hex} for the
expected entropy can be found in
\cite{Wolpert:95,Hutter:01xentropy}. The only thing to be shown is
that $h$ is concave. This may be done by exploiting special
properties of the digamma function $\psi$ (see
\cite[Chp.6]{Abramowitz:74}).
There are fast implementations of $\psi$ and its derivatives and
exact expressions for integer and half-integer arguments
(see Appendix \ref{secPsi} for details).

\fexample{exHexact}{Exact robust expected entropy}{
To see how the derived formulas can be used,
let us compute the upper and lower expected entropy for
for
\beqn\textstyle
   d=2,\quad n_1=3,\quad n_2=6, \qmbox{i.e.} n=9,
   \qmbox{and} s=1,\qmbox{hence} \s={1\over 10}
\eeqn
The general correspondence \req{thtrel} becomes
\beqn\textstyle
  \u_1={3+t_1\over 10},\quad \u_2={6+t_2\over 10},
  \qmbox{hence} \vt^0=\v 0 \qmbox{implies}  \vu^0=\left({0.3\atop 0.6}\right).
\eeqn
Using $n_1<n_2$, \req{Fmin} implies
\beqn\textstyle
  i^{\low H}=2,\quad \vt^{\low H}=\left({0\atop 1}\right),
  \qmbox{hence} \vu^{\low H}=\left({0.3\atop 0.7}\right).
\eeqn
From \req{Fmax}, using $i_1=1$ and $i_2=2$, we get
\beqn\textstyle
  \tilde\u=\min\left\{{1+3\over 9+1},{1+3+6\over 2\cdot(9+1)}\right\}={4\over 10},
  \qmbox{hence} \vu^{\up H}=\max\{\vu^0,\tilde\u\}=\left({0.4\atop 0.6}\right).
\eeqn
This shows that the upper bound is assumed in a/the corner
$\vt^{\up H}=\left({1\atop 0}\right)$.
Inserting these $\u$ into \req{hex}, we get
\beqn\textstyle
h({3\over 10}) = {2761\over 8400},\qquad
h({4\over 10}) = {2131\over 6300},\qquad
h({6\over 10}) = {1207\over 4200},\qquad
h({7\over 10}) = {847\over 3600}.
\eeqn
Putting everything together we get the robust $H$ estimate
\bqan
  \up{\low H} &=& \textstyle [H(\vu^{\low H}),H(\vu^{\up H})]
  = [h({3\over 10})+h({7\over 10})\, ,\, h({4\over 10})+h({6\over 10})]
\\
  &=& \textstyle [{7106\over 12600}\, ,\, {7883\over 12600}]
  \doteq [0.5639,0.6256]
\eqan
The size of this interval is ${37\over 600}$, so $\up H-\low H\doteq
0.0616$ is of the order of $\s$.
}

In general, in order to apply Theorem \ref{thEEI}, we need to be
able to (a) somehow compute $F(\vu)$, e.g.\ compute the expectation
$E_\vt[\cal F]$, (b) verify whether $F(\vu)$ has the form $\sum_i
f(\u_i)$, which is often trivial, e.g.\ if $\cal
F(\vpi)=\sum_i\phi(\pi_i)$, and (c) prove concavity or convexity of
$F$. In the following sections we derive conservative approximations
for more general $F(\vu)$.

\section{Approximate Robust Intervals}\label{secAEI}

In this section we derive approximations for $\up{\low F}$
suitable for arbitrary, twice differentiable functions
$F(\vu)$.
The derived approximations for $\up{\low F}$ will be robust in the
sense of covering set $\up{\low F}$ (for any $n$), and the
approximations will be ``good'' if $n$ is not too small. We do this
by means of a finite Taylor series expansion in $\sigma:={s\over
n+s}$ and by bounding the remainder.

In the following, we treat $\sigma$ as a (small) expansion
parameter. For $\vu,\vu^*\in\Delta'$ we have
\beq\label{dtbnd}
  \u_i-\u_i^* \;=\; \sigma\!\cdot\!(t_i-t_i^*) \qmbox{and}
  |\u_i-\u_i^*| \;=\; \sigma|t_i-t_i^*| \;\leq\; \sigma
  \qmbox{with} \sigma:=\textstyle{s\over n+s}.
\eeq
Hence we may Taylor-expand $F(\vu)$ around $\vu^*$, which leads to a
Taylor series in $\sigma$. This shows that $F$ is approximately
linear in $\vu$ and hence in $\vt$. A linear function on a simplex
assumes its extreme values at the vertices of the simplex. This has
already been encountered in Section \ref{secEEI}. The consideration
above is a simple explanation for this fact. This also shows that
the robust interval $\up{\low F}$ is of size $\up F-\low
F=O(\sigma)$.\footnote{$f(\v n,\vt,s)=O(\sigma^k)$
$\;:\Leftrightarrow\;$ $\exists c\,\forall\v
n\in\SetN_0^d,\,\vt\in\Delta,\,s>0$ : $|f(\v n,\vt,s)|\leq
c\sigma^k$, where $\sigma={s\over n+s}$.} Any approximation to
$\up{\low F}$ should hence be at least $O(\sigma^2)$. The expansion
of $F$ to $O(\sigma)$ is
\beq\label{Fexpand}
  F(\vu) \;=\; \overbrace{F(\vu^*)}^{F_0=O(1)} +
  \overbrace{\sum_i[\partial_i F(\v{\check\u})](\u_i-\u_i^*)}^{F_R=O(\sigma)},
\eeq
where $\partial_i F(\v{\check\u})$ is the partial derivative
$\partial F(\v{\check\u})/\partial\check\u_i$ of
$F(\v{\check\u})$ w.r.t.\ $\check\u_i$. For suitable
$\v{\check\u}=\v{\check\u}(\vu,\vu^*)\in\Delta'$ this expansion is
exact ($F_R$ is the exact remainder).
Natural points for expansion are $t_i^*={1\over d}$ in the center
of $\Delta$, or possibly also $t_i^*={n_i\over n}=\u_i^*$. Here,
we expand around the improper point $t_i^*:=t_i^0\equiv 0$, which
is outside(!) $\Delta$, since this makes expressions particularly
simple.$\!\!$\footnote{The order of accuracy
$O(\sigma^2)$ we will encounter is the same for all choices of $\vu^*$.
The concrete numerical errors differ of course. The
choice $\vt^*=\v 0$ can lead to $O(d)$ smaller $F_R$ than the
natural center point $\vt^*={1\over d}$, but is more likely a
factor $O(1)$ larger. The exact numerical values depend on the
structure of $F$.} Eq.\req{dtbnd} is still valid in this case,
and $F_R$ is exact for some $\v{\check\u}$ in
\beqn
  \Deltapl:=\{\vu\,:\,\u_i\geq\u_i^0\,\forall i,\;\u_\p\leq 1\},
  \qmbox{where}
  \u_i^0={n_i\over n+s}.
\eeqn
Note that we keep the exact condition
$\vu\in\Delta'$. $F$ is usually already defined on
$\Deltapl$ or extends from $\Delta'$ to $\Deltapl$
without effort in a natural way (analytical continuation). We
introduce the notation
\beq\label{eqleqsq}
  F\leqsq G \qquad:\Leftrightarrow\qquad
  F\leq G \qmbox{and} F=G+O(\sigma^2),
\eeq
stating that $G$ is a ``good'' upper bound on $F$.
The following bounds hold for arbitrary differentiable functions.
In order for the bounds to be ``good,$\!$'' $F$ has to be Lipschitz
differentiable in the sense that there exists a constant $c$ such that
\beqn
  |\partial_i F(\vu)|\leq c \qmbox{and}
  |\partial_i F(\vu)-\partial_i F(\vu')|
  \leq c|\vu-\vu'|
\eeqn
\beq\label{Lipschitz}
  \forall\,\vu,\vu'\in\Deltapl
  \qmbox{and} \forall\,i\in\{1,...,d\}
\eeq
\beqn
  \mbox{If $F$ depends also on $\v n$, e.g.\ via $\sigma$ or
  $\vu^0$, then $c$ shall be independent of them.}
\eeqn
The Lipschitz condition is satisfied, for instance, if the
curvature $\partial^2 F$ is uniformly bounded. This is satisfied
for the expected entropy $H$ (see \req{hex}), but violated for
the approximation $E_\vt[{\cal H}]\approx{\cal H}(\vu)$ if $n_i=0$ for
some $i$.

\begin{theorem}[Approximate robust intervals]\label{thARI}
Assume $F:\Deltapl\to\SetR$ is a Lipschitz differentiable
function \req{Lipschitz}.
Let $[\low F,\up F]$ be the global [minimum,maximum] of $F$ restricted
to $\Delta'$. Then
\beqn
  F(\vu^1) \;\leqsq\; \up F \;\leqsq\; F_0+F_R^{ub}
  \;\;\mbox{where}\;\; F_R^{ub}:=\max_i F_{iR}^{ub}
  \;\mbox{and}\;\;
  F_{iR}^{ub} \;:=\; \sigma\maxo_{\vu\in\Deltapl}
    [\partial_i F(\vu)],
\eeqn\vspace{-1ex}
\beqn
  F_0+F_R^{lb} \;\leqsq\; \low F \;\leqsq\; F(\vu^2)
  \;\;\mbox{where}\;\; F_R^{lb}:=\min_i F_{iR}^{lb}
  \;\;\mbox{and}\;\;
  F_{iR}^{lb} \;:=\; \sigma\mino_{\vu\in\Deltapl}
    [\partial_i F(\vu)],
\eeqn
$F_0:=F(\vu^0)$, and
$t^1_i:=\delta_{ii^1}$ with $i^1:=\arg\max_i F_{iR}^{ub}$, and
$t^2_i:=\delta_{ii^2}$ with $i^2:=\arg\min_i F_{iR}^{lb}$, and $\leqsq$
defined in \req{eqleqsq} means $\leq$ {\em and} $=\;+O(\sigma^2)$,
where $\sigma=1-\u_\p^0$.
\end{theorem}
For conservative estimates, the lower bound on $\low F$ and the
upper bound on $\up F$ are the interesting ones. Together with the
``inner'' bounds $F(\vu^1)$ and $F(\vu^2)$, they also yield
interesting information about the accuracy of the approximations:
$F_0+F_R^{ub}-F(\vu^1)$ is an upper bound on the (unknown)
approximation error $F_0+F_R^{ub}-\up F$, and similarly for $\low
F$.

\paradot{Proof} We start by giving an $O(\sigma^2)$ bound on
$\up F_R=\maxo_{\vu\in\Delta'}F_R(\vu)$. We first insert
\req{dtbnd} with $\vt^*=\vt^0\equiv\v 0$ into
\req{Fexpand} and treat $\v{\check\u}$ and $\vt$ as separate
variables:
\beqn
  F_R(\v{\check\u},\vt) \;=\; \sigma\sum_i
  [\partial_i F(\v{\check\u})]\cdot t_i
  \;\leqsq\; \maxo_{\v{\check\u}\in\Deltapl}
  \bigg\{\sigma\sum_i[\partial_i F(\v{\check\u})]\cdot t_i\bigg\}
  \;\leqsq\;
  \sum_i F_{iR}^{ub}\cdot t_i
\eeqn
\beq\label{defFiRub}
  \qmbox{with}
  F_{iR}^{ub} \;:=\; \sigma\maxo_{\v{\check\u}\in\Deltapl}
    [\partial_i F(\v{\check\u})]
\eeq
The first inequality is obvious, the second follows from the
convexity of $\max$. From assumption \req{Lipschitz} we get
$\partial_i F(\vu)-\partial_i F(\vu') = O(\sigma)$ for all
$\vu,\vu'\in\Deltapl$, since $\Deltapl$ has
dia\-meter $O(\sigma)$. Due to one additional $\sigma$ in
\req{defFiRub} the expressions in \req{defFiRub} change only by
$O(\sigma^2)$ when introducing or dropping $\maxo_{\v{\check\u}}$
anywhere. This shows that the inequalities are tight within
$O(\sigma^2)$ and justifies $\leqsq$. We now upper bound $F_R(\vu)$:
\beq\label{FRup}
  \up F_R = \maxo_{\vu\in\Delta'} F_R(\vu)
  \leqsq \maxo_{\vt\in\Delta}\maxo_{\v{\check\u}\in\Deltapl}
  F_R(\v{\check\u},\vt)
  \leqsq
  \maxo_{\vt\in\Delta}\sum_i F_{iR}^{ub}\cdot t_i
  =
  \max_i F_{iR}^{ub} =: F_R^{ub}
\eeq
A linear function on $\Delta$ is maximized by setting the $t_i$
component with largest coefficient to 1. This shows the last
equality. The maximization over $\v{\check\u}$ in \req{defFiRub} can
often be performed analytically, leaving an easy $O(d)$ time task
for maximizing over $i$.

We have derived an upper bound $F_R^{ub}$ on $\up F_R$. Let us
define the corner $t_i=\delta_{ii^1}$ of $\Delta$ with
$i^1:=\arg\max_i F_{iR}^{ub}$. Since $\up F_R\geq F_R(\vu)$
for all $\vu$, $F_R(\vu^1)$ in particular is a lower bound
on $\up F_R$. A similar line of reasoning as above shows that that
$F_R(\vu^1)=\up F_R+O(\sigma^2)$. Using $\up{F+const.}=\up
F+const.$ we get $O(\sigma^2)$ lower and upper bounds on $\up F$,
i.e.\ $F(\vu^1)\leqsq\up F\leqsq F_0+F_R^{ub}$. $\low F$ is bound
similarly with all max's replaced by min's and inequalities
reversed. Together this proves the Theorem \ref{thARI}. \qed

In the following sections we assume the definitions/notation of
Theorem \ref{thARI} for $F$ and analogous ones for all other
occurring estimators ($G,H,I,...$).

\section{Error Propagation}\label{secEP}

We now show how bounds of elementary functions obtained by Theorem
\ref{thARI} can be used to get bounds for more complex composite
functions, especially for sums and products of functions. The
results are used in Section \ref{secIEMI} for deriving robust
intervals for the mutual information for which exact solutions are
not known.

\paradot{Approximation of $\up{\low F}$ (special cases)}
For the special case $F(\vu)=\sum_i f(\u_i)$ we have $\partial_i
F(\vu)=f'(\u_i)$. For concave $f$ like in case of the entropy we
get particularly simple bounds
\bqa\label{defFRulb}
  \nq F_{iR}^{ub} &= \sigma\maxo_{\vu\in\Deltapl} f'(\u_i)
   = \sigma f'(\u_i^0),\quad\;\;
  & F_R^{ub} = \sigma\max_i f'(\u_i^0)
  = \sigma f'(\textstyle{\min_i n_i\over n+s}),\quad\;\;
\\ \nonumber
  \nq F_{iR}^{lb} &= \sigma\mino_{\vu\in\Deltapl} f'(\u_i)
  = \sigma f'(\u_i^0+\sigma),\quad
  & F_R^{lb} = \sigma\min_i f'(\u_i^0+\sigma)
  = \sigma f'(\textstyle{\max_i n_i+s\over n+s}),
\eqa
where we have used $\maxo_{\vu\in\Deltapl} f'(\u_i)
=\max_{\u_i\in[\u_i^0,\u_i^0+\sigma]} f'(\u_i)=f'(\u_i^0)$, and similarly
for $\min$. Analogous results hold for convex functions. In case
the maximum cannot be found exactly one is allowed to further
increase $\Deltapl$ as long as its diameter remains $O(\sigma)$.
Often an increase to
$\Box':=\{\vu:\u_i^0\leq\u_i\leq\u_i^0+\sigma\} \supset
\Deltapl \supset \Delta'$ makes the problem easy. Note that if we
were to perform these kind of crude enlargements on
$\maxo_{\vu}F(\vu)$ directly we would loose the bounds by
$O(\sigma)$.

\fexample{exHbounds}{Approximate robust expected entropy}{
Let us compare the exact robust estimate of the expected entropy for
$n_1=3$, $n_2=6$, $s=1$ (hence $n=9$, and $\sigma={1\over 10}$)
computed in Example \ref{exHexact} with this approximation:
Using the expressions for $h'$ from Appendix
\ref{secPsi}, we get
\beqn\textstyle
  h'({3\over 10})= {13051\over 2520}-\fr12\pin^2 \qmbox{and}
  h'({7\over 10})= {91717\over 8400}-{7\over 6}\pin^2,
\eeqn
where $\pin\doteq 3.1415$. From \req{corEEI} and \req{defFRulb}
we get
\beqn\textstyle
  H_0=H(\vu^0)=h({3\over 10})+h({6\over 10})={69\over 112},\qquad
  H_R^{ub}={1\over 10}h'({3\over 10}),\qquad
  H_R^{lb}={1\over 10}h'({7\over 10}).
\eeqn
Together with the expressions from Example \ref{exHexact} we get
the conservative estimate
\beqn
  [H_0+H_R^{lb}\, ,\, H_0+H_R^{ub}] \doteq [0.5564,0.6404].
\eeqn
The approximation accuracy
\beqn
  H_0+H_R^{ub}-\up H \doteq 0.0148 \qmbox{and}
  \low H-H_0-H_R^{lb} \doteq 0.0074
\eeqn
is consistent with our $O(\s^2)$ estimation. If exact expressions
are not available we can upper bound the widening by
\beqn
  H_0+H_R^{ub}-H(\vu^1) \doteq 0.0148 \qmbox{and}
  H(\vu^2)-H_0-H_R^{lb} \doteq 0.0074
\eeqn
Since generally $\vu^2=\vu^{\low H}$ and in our example also
$\vu^1=\vu^{\up H}$, the numbers coincide.
}

\begin{figure}
\includegraphics[width=0.5\textwidth]{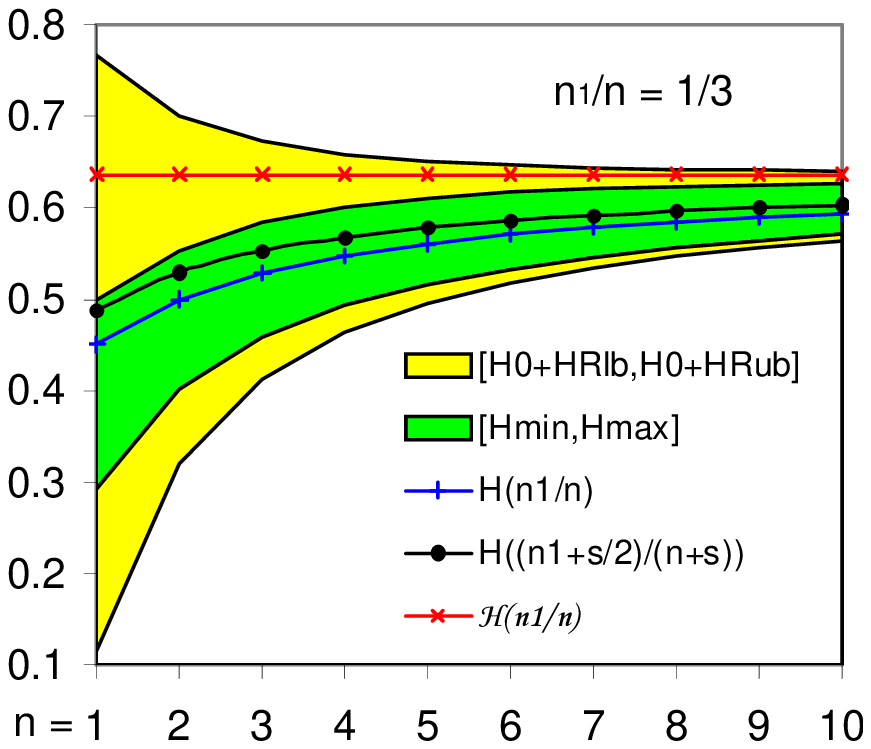}
\includegraphics[width=0.5\textwidth]{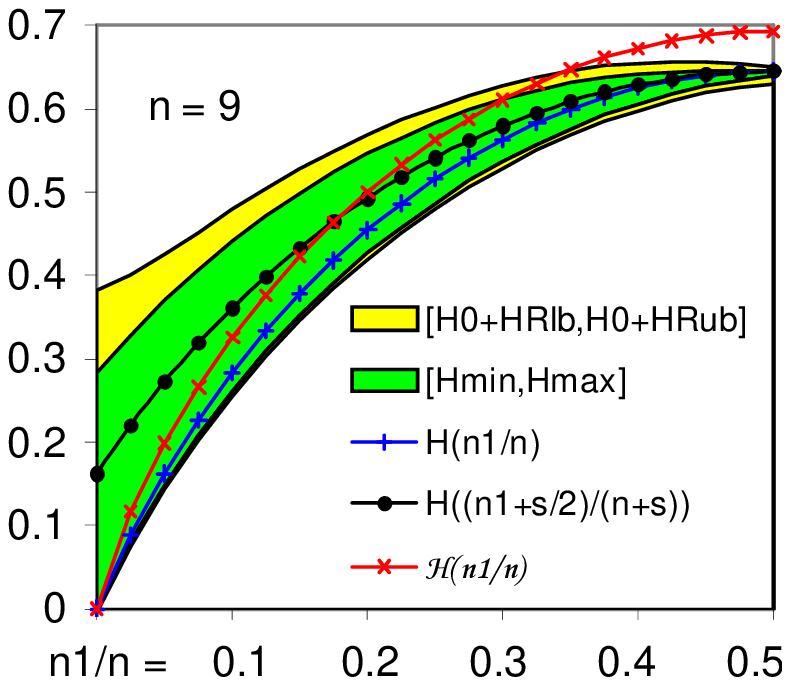}
\caption{\label{figHnx}[Expected Entropy] The figures display the
various (expected) entropy estimates for $s=1$: The left figure for
$n_1/n=1/3$ and $n=1...10$. The right figure for $n=9$ and
$n_1/n=0...0.5$. The ``intersection'' $n_1=3$ and $n_2=6$ is treated
analytically in Examples \ref{exHexact} and \ref{exHbounds}. The
green (dark gray) area is the exact robust interval $[\low H\, ,\,
\up H]$ from Theorem \ref{corEEI}. The yellow+green (gray) area is the
conservative estimate $[H_0+H_R^{lb}\, ,\, H_0+H_R^{ub}]$ from
Theorem \ref{thARI}. The area $[H(\v u^2)\, ,\,H(\v u^1)]$ is not
shown, since (here) it essentially coincides with $\up{\low H}$).
Some point estimates $H({\v n\over n})$, $H({\v n+\v 1/2\over
n+1})$, and ${\cal H}({\v n\over n})$ are also shown.}
\end{figure}

\fexample{exHfig}{Entropy: dependency on $\v n$}{
Figure \ref{figHnx} (left) shows how the size of the (conservative)
robust interval of the expected entropy $H$ varies with the sample
size $n$. We considered $s=1$ and $d=2$ and kept $n_1/n=\frs13$ and
$n_2/n=\frs23$ fix (allowing for fractional $\v n$). We clearly see
that the yellow (light gray) region diminishes quickly compared to
the green (dark gray) region with increasing
$n$, i.e.\ the approximation accuracy gets better for larger $n$.
Some point estimates $H({\v n\over n})$, $H({\v n+\v 1/2\over
n+1})$, and ${\cal H}({\v n\over n})$ are also shown.
Figure \ref{figHnx} (right) shows the intervals for fixed $n=9$,
while varying $n_1/n=0...0.5$  ($n_1/n=0.5...1$ is symmetric).
The interval $\up{\low H}$ is shorter for more uniform $\vu$,
since $H$ (like $\cal H$) varies more closer to the boundary of
$\Delta$.
The $[H(\v u^2),H(\v u^1)]$ region is not shown since it is
identical to $\up{\low H}$ (also in the left graph except for
$n=1$).
For $n=9$ and $n_1/n=1/3$ we recover the results of Examples
\ref{exHexact} and \ref{exHbounds} (left and  right figure). }

\paradot{Error propagation}
Assume we found bounds for estimators $G(\vu)$ and
$H(\vu)$ and we want now to bound the sum
$F(\vu):=G(\vu)+H(\vu)$. In the direct approach $\up
F\leq \up G+\up H$ we may lose $O(\sigma)$. A simple example is
$G(\vu)=\u_i$ and $H(\vu)=-\u_i$ for which
$F(\vu)=0$, hence $0=\up F\leq \up G+\up
H=\u_i^0+\sigma-\u_i^0=\sigma$, i.e.\ $\up F\not\leqsq \up G+\up
H$.
We can exploit the techniques of the previous section to obtain
$O(\sigma^2)$ approximations.
\beqn
  F_{iR}^{ub} \;=\;
  \sigma\maxo_{\vu\in\Deltapl}\partial_i F(\vu)
  \;\leqsq\;
  \sigma\maxo_{\vu\in\Deltapl}\partial_i G(\vu) +
  \sigma\maxo_{\vu\in\Deltapl}\partial_i H(\vu)
  \;=\; G_{iR}^{ub} + H_{iR}^{ub}
\eeqn

\begin{theorem}[Error propagation: Sum]\label{thEps}
Let $G(\vu)$ and $H(\vu)$ be Lipschitz differentiable and
$F(\vu)=\alpha G(\vu)+\beta H(\vu)$, $\alpha,\beta\geq 0$, then $\up
F\leqsq F_0+F_R^{ub}$ and $\low F\geqsq F_0+F_R^{lb}$, where
$F_0=\alpha G_0+\beta H_0$, and $F_{iR}^{ub}\leqsq \alpha
G_{iR}^{ub}+\beta H_{iR}^{ub}$, and $F_{iR}^{lb}\geqsq \alpha
G_{iR}^{lb}+\beta H_{iR}^{lb}$.
\end{theorem}
It is important to notice that $F_R^{ub}\not\leqsq
G_R^{ub}+H_R^{ub}$ (use previous example),
i.e.\ $\max_i[G_{iR}^{ub}+H_{iR}^{ub}]\not\leqsq\max_i
G_{iR}^{ub}+\max_i H_{iR}^{ub}$. $\max_i$ can not be pulled in and
it is important to propagate $F_{iR}^{ub}$, rather than
$F_R^{ub}$.

Every function $F$ with bounded curvature can be written as a sum of
a concave function $G$ and a convex function $H$. For convex and
concave functions, determining bounds is particularly easy, as we
have seen. Often $F$ decomposes naturally into convex and concave
parts as is the case for the mutual information, addressed later.
Bounds can also be derived for products.
%
\begin{theorem}[Error propagation: Product]\label{thEpp}
Let $G,H:\Deltapl\to[0,\infty)$ be non-nega\-tive
Lipschitz differentiable
functions \req{Lipschitz} with non-negative
derivatives $\partial_i G,\partial_i H\geq 0$ $\forall i$ and
$F(\vu)=G(\vu)\cdot H(\vu)$, then $\up F\leqsq
F_0+F_R^{ub}$, where $F_0= G_0\cdot H_0$, and $F_{iR}^{ub}\leqsq
G_{iR}^{ub}(H_0+H_R^{ub})+ (G_0+G_R^{ub})H_{iR}^{ub}$, and
similarly for $\low F$.
\end{theorem}

\paradot{Proof} We have
\bqan
  & \nq\nq F_{iR}^{ub}
  \;=\; \sigma\max\partial_i F
  \;=\; \sigma\max\partial_i (G\!\cdot\!H)
  \;=\; \sigma\max[(\partial_i G)H+G(\partial_i H)]
  \;\leqsq\;
\\
  & \nq\nq \sigma(\max\partial_i G)(\max H) +
  \sigma(\max G)(\max \partial_i H)
  \leqsq G_{iR}^{ub}(H_0\!+\!H_R^{ub}) + (G_0\!+\!G_R^{ub})H_{iR}^{ub}
\eqan
where all functions depend on $\vu$ and all $\max$ are over
$\vu\in\Deltapl$. There is one subtlety in the last
inequality: $\max G\neq \up G\leqsq G_0+G_R^{ub}$.
The reason for the $\neq$ being that the maximization is taken
over $\Deltapl$, not over $\Delta'$ as in the definition
of $\up G$. The
correct line of reasoning is as follows:
\beqn
  \maxo_{\vu\in\Deltapl} G_R(\vu)
  \leqsq
  \maxo_{\vt\in\Deltal}\sum_i G_{iR}^{ub}\cdot t_i
  =
  \max\{0,\max_i G_{iR}^{ub}\} = G_R^{ub}
  \;\Rightarrow\;
  \max G\leqsq G_0+G_R^{ub}
\eeqn
The first inequality can be proven in the same way as \req{FRup}. In
the first equality we set the $t_i=1$ with maximal $G_{iR}^{ub}$
{\em if} it is positive. If all $G_{iR}^{ub}$ are negative we set
$\vt\equiv\v 0$. We assumed $G\geq 0$ and $\partial_i G\geq 0$,
which implies $G_R\geq 0$. So, since $G_R\geq 0$ anyway, this
subtlety is ineffective. Similarly for $\max H_R$.\qed

It is possible to remove the rather strong non-negativity
assumptions. Propagation of errors for other combinations like
ratios $F=G/H$ may also be obtained.

\section{Robust Intervals for Expected Mutual Information}\label{secIEMI}

We illustrate the application of the previous results on the Mutual
Information between two random variables $\imath\in\{1,...,d_1\}$
and $\jmath\in \{1,...,d_2\}$.

\paradot{Mutual Information}
Consider an i.i.d.\ random process with outcome
$(i,j)\in\{1,...,d_1\}\times\{1,...,d_2\}$ having joint probability
$\pi_{ij}$, where $\vpi\in\Deltapi :=\{\v x\in\SetR^{d_1\times
d_2}\,:\,x_{ij}\geq 0\,\forall ij,\; x_\pp=1\}$. An important
measure of the stochastic dependence of $\imath$ and $\jmath$ is the
mutual information
\bqa\nonumber
  {\cal I}(\vpi) \;=\; \sum_{i=1}^{d_1}\sum_{j=1}^{d_2}
  \pi_{ij}\log{\pi_{ij}\over\pi_{i\p}\pi_{\p j}}
  &=& \sum_{ij}\pi_{ij}\log\pi_{ij} -\!
      \sum_{i}\pi_{i\p}\log\pi_{i\p} -\!
      \sum_{j}\pi_{\p j}\log\pi_{\p j}
\\ \label{mi}
  &=& {\cal H}(\vpi_{\imath\p}) +
      {\cal H}(\vpi_{\p\jmath}) - {\cal
      H}(\vpi_{\imath\jmath}),
\eqa
where $\pi_{i\p}=\sum_j\pi_{ij}$ and $\pi_{\p j}=\sum_i\pi_{ij}$ are row
and column marginal chances. Again, we assume a Dirichlet prior
over $\vpi_{\imath\jmath}$, which leads to a Dirichlet posterior
$p(\vpi_{\imath\jmath}|\v
n)\propto\prod_{ij}\pi_{ij}^{n_{ij}+st_{ij}-1}$ with $\vt\in
\Delta$. The expected value of $\pi_{ij}$ is
\beqn
  E_\vt[\pi_{ij}]={n_{ij}+st_{ij}\over n+s}=:\u_{ij}
\eeqn
The marginals $\vpi_{i\p}$ and $\vpi_{\p j}$ are also Dirichlet with
expectation $\u_{i\p}$ and $\u_{\p j}$.
The expected mutual information $I(\vu):=E_\vt[{\cal I}]$
can, hence, be expressed in terms of the expectations of
three entropies $H(\vu):=E_\vt[{\cal H}]$ (see \req{hex})
\bqan
  I(\vu) &=& H(\vu_{\imath\p})+H(\vu_{\p\jmath})-H(\vu_{\imath\jmath})
 \;=\; H_{row}+H_{col}-H_{joint}
\\
  &=& \sum_i h(\u_{i\p})+\sum_j h(\u_{\p j})-\sum_{ij}h(\u_{ij}),
\eqan
where here and in the following we index quantities with $joint$,
$row$, and $col$ to denote to which distribution the quantity
refers.

\paradot{Crude bounds for $I(\vu)$}
Estimates for the robust IDM interval
$[\mino_{\vt\in\Delta}E_\vt[{\cal I}] \,,\, \maxo_{\vt\in\Delta
}E_\vt[{\cal I}]]$ can be obtained by [minimizing,maximizing]
$I(\vu)$. A crude upper bound can be obtained as
\bqan
  \up I &:=& \maxo_{\vt\in\Delta} I(\vu) \;=\;
  \max[H_{row}+H_{col}-H_{joint}] \;\leq\;
\\
  & & \max H_{row} + \max H_{col} - \min H_{joint} \;=\;
  \up H_{row} + \up H_{col} -  \low H_{joint},
\eqan
where exact solutions to $\up H_{row}$, $\up H_{col}$ and
$\low H_{joint}$ are available from Section \ref{secEEI}.
Similarly
$
  \low I \geq \low H_{row} + \low H_{col} -  \up H_{joint}
$. The problem with these bounds is that, although good in some
cases, they can become arbitrarily crude.
The following $O(\sigma^2)$ bound can be derived by exploiting the
error sum propagation Theorem \ref{thEps}.

\begin{theorem}[Bound on lower and upper expected Mutual Information]\label{thMIbnd}
The following bounds on the expected mutual information
$I(\vu)=E_\vt[{\cal I}]$ are valid:
\bqan
 & & \nq I(\vu^1) \leqsq \up I\leqsq I_0+I_R^{ub} \qmbox{and}
 I_0+I_R^{lb} \leqsq \low I \leqsq  I(\vu^2),
 \qmbox{where} \\
 & & \nq I_0=I(\vu^0)=H_{0row}+H_{0col}-H_{0joint} =
 \textstyle \sum_i h(\u_{i\p}^0)+ \sum_j h(\u_{\p j}^0)- \sum_{ij} h(\u_{ij}^0), \\
 & & \nq I_{ijR}^{ub}\leqsq
 H_{iRrow}^{ub}+H_{jRcol}^{ub}-H_{ijRjoint}^{lb} =
 h'(\u_{i\p}^0)+h'(\u_{\p j}^0)-h'(\u_{ij}^0\!+\!\sigma), \\
 & & \nq I_{ijR}^{lb}\geqsq
 H_{iRrow}^{lb}+H_{jRcol}^{lb}-H_{ijRjoint}^{ub} =
 h'(\u_{i\p}^0\!+\!\sigma)+h'(\u_{\p j}^0\!+\!\sigma)-h'(\u_{ij}^0),
\eqan
with $h$ defined in \req{hex}, and $t^0_{ij}=0$, and
$t^1_{ij}=\delta_{(ij)(ij)^1}$ with $(ij)^1=\arg\max_{ij}
I_{ijR}^{ub}$, and $t^2_{ij}=\delta_{(ij)(ij)^2}$ with
$(ij)^2=\arg\min_{ij} I_{ijR}^{lb}$, and $I_R^{ub}=\max_{ij}I_{ijR}^{ub}$,
and $I_R^{lb}=\max_{ij}I_{ijR}^{lb}$.
\end{theorem}

\section{The IDM for Product Spaces}\label{secPS}

In the last section we considered the ``full'' IDM on the product of
two random variables. The structure of the problem suggests
considering a smaller ``product'' of IDMs as described below, which
can lead to better estimates.

Product spaces $\Omega=\Omega_1\times...\times\Omega_m$ with
$\Omega_k=\{1,...d_k\}$ occur frequently in practical problems,
e.g.\ in the mutual information ($m=2$), in robust trees ($m=3$),
or in Bayesian nets in general ($m$ large). Without loss of
generality we only discuss the $m=2$ case in the following.
Ignoring the underlying structure in $\Omega$, a Dirichlet prior
in case of unknown chances $\pi_{\imath\jmath}$ and an IDM as used
in Section \ref{secIEMI} with
\beq\label{IDMfull}
  \vt\in\Delta :=\{\vt\in\SetR^{d_1\times
  d_2}\equiv\SetR^{d_1}\otimes\SetR^{d_2}\,:\,t_{ij}\geq
  0\;\forall\,ij,\; t_\pp=1\}
\eeq
seems natural.

On the other hand, if we take into account the structure of $\Omega$
and go back to the original motivation of the IDM, this choice is
far less obvious. Recall that one of the major motivations of the
IDM was its representation invariance in the sense that
inferences are not affected when grouping or splitting events in
$\Omega$. For unstructured spaces like $\Omega_k$ this is a
reasonable principle. For illustration, let us consider objects of
various {\em shape} and {\em color}, i.e.\
$\Omega=\Omega_1\times\Omega_2$, $\Omega_1=\{ball, pen, die, ...\}$,
$\Omega_2=\{yellow, red, green, ...\}$ in generalization to Walley's
bag of marbles example. Assume we want to detect a potential
dependency between {\em shape} and {\em color} by means of their
mutual information $I$. If we have no prior idea on the possible
kind of colors, a model which is independent of the choice of
$\Omega_2$ is welcome. Grouping red and green, for instance,
corresponds to grouping $(x_{i1}$, $x_{i2}$, $x_{i3}$, $x_{i4},
...)$ to $(x_{i1}$, $x_{i2}+x_{i3}$, $x_{i4}, ...)$ {\em for all
shapes $i$}, where $\v x\in\{\v n,\vpi,\vt,\vu\}$. Similarly for the
different shapes, for instance we could group all round or all
angular objects. The ``smallest IDM'' which respects this invariance
is the one which considers all
\beq\label{IDMprod}
  \vt\in\DeltaOX:=\Delta_{d_1}\otimes\Delta_{d_2}
  \;\propersubset\;\Delta.
\eeq
The tensor or outer product $\otimes$ is defined as $(\v v\otimes\v
w)_{ij}:=v_iw_j$ and $V\otimes W:=\{\v v\otimes\v w:\v v\in V,\, \v
w\in W\}$. It is a bilinear (not linear!) mapping. This smaller
product IDM $\DeltaOX$ is invariant under arbitrary grouping of
columns and rows of the chance matrix $(\vpi_{ij})_{1\leq i\leq
d_1,1\leq j\leq d_2}$. In contrast to the larger full IDM $\Delta$
it is not invariant under arbitrary grouping of matrix cells, but
there is anyway little motivation for the necessity of such a
general invariance. General non-column/row cross groupings would
destroy the product structure of $\Omega$ and with that the mere
concepts of shape and color, and their correlation. For $m>2$ as in
Bayes-nets cross groupings look even less natural. Whether the
$\DeltaOX$ or the larger simplex $\Delta$ is the more appropriate
IDM depends on whether one regards the structure
$\Omega_1\times\Omega_2$ of $\Omega$ as a natural prior knowledge or
as an arbitrary a posteriori choice. The smaller IDM has the
potential advantage of leading to more precise predictions (smaller
robust sets).

Let us consider an estimator $F:\Delta\to\SetR$ and its
restriction $F_\ots:\DeltaOX\to\SetR$.
Robust intervals $[\low F,\up F]$ for $\Delta$ are generally wider
than robust intervals $[\low{F}_\ots,\up F_\ots]$ for $\DeltaOX$.
Fortunately not much. Although $\DeltaOX$ is a {\em
lower-dimensional} subspace of $\Delta$, it contains all vertices of
$\Delta$. This is possible since $\DeltaOX$ is a {\em nonlinear}
subspace. The set of ``vertices'' in both cases is $\{\vt\,:\,t_{ij}
= \delta_{ii_0}\delta_{jj_0},\;i_0\in\Omega_1,\;j_0\in\Omega_2\}$.
Hence, {\em if} the robust interval boundaries $\up{\low F}$ are
assumed in the vertices of $\Delta$ {\em then} the
interval for the $\DeltaOX$ IDM model is the
same ($\up{\low F}=\up{\low F}_\ots$). Since the condition is
``approximately'' true, the conclusion is
``approximately'' true. More precisely:

\begin{theorem}[IDM bounds for product spaces]
The $O(\sigma^2)$ bounds of Theorem \ref{thARI} on the robust
interval $\up{\low F}$ in the full IDM $\Delta$ \req{IDMfull},
remain valid for $\up{\low F}_\ots$ in the product IDM $\DeltaOX$
\req{IDMprod}.
\end{theorem}

\paradot{Proof}
\beqn
  F(\vu^1) \leq \up F_\ots \leq \up F \leq
  F_0+F_R^{ub} = F(\vu^1)+O(\sigma^2),
\eeqn
where $\up F_\ots:=\maxo_{\vt\in\deltaOX}F(\vu)$ and $\vu^1$
was the ``$F_R$ maximizing'' vertex as defined in Theorem
\ref{thMIbnd} ($F(\vu^1)\leqsq\up F$). The first inequality
follows from the fact that all $\Delta$ vertices also belong to
$\DeltaOX$, i.e.\ $\vt^1\in\DeltaOX$. The second inequality
follows from $\DeltaOX\subset\Delta$. The remaining (in)equalities
follow from Theorem \ref{thARI}. This shows that $|\up F_\ots-\up
F|=O(\sigma^2)$, hence $F_0+F_R^{ub}$ is also an $O(\sigma^2)$
upper bound to $\up F_\ots$. This implies that to the
approximation accuracy we can achieve, the choice between $\Delta$
and $\DeltaOX$ is irrelevant. \qed

\section{Robust Credible Intervals}\label{secCI}

So far we have considered robust intervals of {\em expected} values
$F=E_\vt[{\cal F}]$. We now briefly consider the problem of how to
combine Bayesian credible intervals for $\cal F$ with robust
intervals of the IDM.

\paradot{Bayesian credible sets/intervals}
For a probability density $p:\SetR^d\to[0,1]$, an $\alpha$-credible
region is a measurable set $A$ for which $p(A):=\int
p(x)\indfct_A(x) d^dx\geq\alpha$, where $\indfct_A(x)=1$ if $x\in A$
and $0$ otherwise, i.e.\ $x\in A$ with probability at least
$\alpha$. For given $\alpha$, there are many choices for $A$. Often
one is interested in ``small'' sets, where the size of $A$ may be
measured by its volume $\Vol(A):=\int\indfct_A(x)d^dx$. Let us
define a/the smallest $\alpha$-credible set
\beqn
  A^{min}:=\mathop{\arg\min}_{A:p(A)\geq\alpha}\Vol(A)
\eeqn
with ties broken arbitrarily. For unimodal $p$, $A^{min}$ can be
chosen as a connected set. For $d=1$ this means that $A^{min}=[a,b]$
with $\int_a^b p(x)dx=\alpha$ is a minimal length highest density
$\alpha$-credible interval. If, additionally $p$ is symmetric around
$E[x]$, then $A^{min}=[E[x]-c,E[x]+c]$ is also symmetric around
$E[x]$.

\paradot{Robust credible sets}
If we have a set of probability distributions $\{p_t(x)$, $t\in
T\}$, we can choose for each $t$ an $\alpha$-credible set $A_t$ with
$p_t(A_t)\geq\alpha$, a minimal one being
$A_t^{min}:=\arg\min_{A:p_t(A)\geq\alpha}\Vol(A)$. A robust
$\alpha$-credible set is a set $A$ which contains $x$ with
$p_t$-probability at least $\alpha$ for {\em all} $t$. A minimal
size robust $\alpha$-credible set is
\beq\label{eqRCI}
  A^{min}:=\mathop{\arg\min}_{A={\textstyle\cup}_t A_t:p_t(A_t)\geq\alpha}
  \Vol(A).
\eeq
It is not easy to deal with this expression, since $A^{min}$ is
{\em not} a function of $\{A_t^{min}:t\in T\}$, and especially does
not coincide with $\bigcup_t A_t^{min}$ as one might expect.

\paradot{Robust credible intervals}
This can most easily be seen for univariate symmetric
unimodal distributions, where $t$ is a translation, e.g.\
$p_t(x)=\mbox{Normal}(E_t[x]=t,\sigma=1)$ with 95\% credible
intervals $A_t^{min}=[t-2,t+2]$. For, e.g.\ $T=[-1,1]$ we get
$\bigcup_t A_t^{min}=[-3,3]$. The credible intervals {\em move}
with $t$. One can get a smaller union if we take the intervals
$A_t^{sym}=[-c_t,c_t]$ symmetric around 0. Since $A_t^{sym}$ is a
non-central interval w.r.t.\ $p_t$ for $t\neq 0$, we have $c_t>2$, i.e.\
$A_t^{sym}$ is larger than $A_t^{min}$, but one can show that the
increase of $c_t$ is smaller than the shift of $A_t^{min}$ by $t$,
hence we save something in the union. The optimal choice is
neither $A_t^{sym}$ nor $A_t^{min}$, but something in-between.

To illustrate this point numerically consider
triangular distributions instead of Gaussians:
\beqn
  p_t(x):=\max\{0\,,\,1\!-\!|x\!-\!t|\},\qquad
  t\in T:=[-\gamma,\gamma],\qquad
  \gamma>0,
\eeqn
\beqn
  \Rightarrow\;
  p_t([a,b])=\left|b^*(1\!-\!\fr12|b^*|)-a^*(1\!-\!\fr12|a^*|)\right|
  \;\mbox{with}\;
  {a^*=\min\{\max\{a,0\},1\}\!-\!t, \atop
  b^*=\min\{\max\{b,0\},1\}\!-\!t.}
\eeqn
One can derive the following expressions for the
$\alpha$-credible intervals, valid for (the interesting case of)
$\alpha\geq\fr12$.

\beqn
  A_t^{min}=[t-1+\sqrt{1-\alpha} \,,\, t+1-\sqrt{1-\alpha}],
\eeqn
\beqn
  \bigcup_{t\in T} A_t^{min}=[-\gamma-1+\sqrt{1-\alpha} \,,\,
  \gamma+1-\sqrt{1-\alpha}].
\eeqn
\beqn
  A^{min}=\left\{
  \begin{array}{ccc}
    \left[-1\!+\!\sqrt{1\!-\!\alpha\!-\!\gamma^2} \,,\,
      1\!-\!\sqrt{1\!-\!\alpha\!-\!\gamma^2}\right] & \mbox{for} &
      \gamma^2\leq\fr12(1\!-\!\alpha), \\
       \left[-\gamma\!-\!1\!+\!\sqrt{2(1\!-\!\alpha)} \,,\,
    \gamma\!+\!1\!-\!\sqrt{2(1\!-\!\alpha)}\right]  &
     \mbox{for} &
      \gamma^2\geq\fr12(1\!-\!\alpha).
  \end{array}\right.
\eeqn
It is easy to see that $A^{min}\subset\bigcup_t A_t^{min}$ and
that $A^{min}$ is a proper subinterval of $\bigcup_t
A_t^{min}$ of shorter length for every $\gamma>0$ and
$\fr12\leq\alpha<1$.

An interesting open question is under which general conditions we
can expect $A^{min}\subseteq\bigcup_t A_t^{min}$. In any case,
$\bigcup_t A_t$ can be used as a conservative estimate for a
robust credible set, since $p_t(\bigcup_{t'} A_{t'})\geq
p_t(A_t) \geq \alpha$ for all $t$.

A special (but important) case which falls outside the above
framework are one-sided credible intervals, where only $A_t$ of
the form $[a,\infty)$ are considered. In this case
$A^{min}=\bigcup_t A_t^{min}$, i.e.\ $A^{min}=[a_{min},\infty)$
with $a_{min}=\max\{a:p_t([a,\infty])\geq\alpha\forall t\}$.

\paradot{Approximations}%
For complex distributions like for the mutual information we have
to approximate \req{eqRCI} somehow. We use the following notation
for shortest $\alpha$-credible {\em intervals} w.r.t.\
a univariate distribution $p_t(x)$:
\beqn
  {\mathop{\widetilde{x}}\limits_\sim}\!\,_t \;\equiv\;
  [ {\mathop{x}\limits_\sim}\!\,_t , \widetilde{x}_t] \;\equiv\;
  [ E_t[x]-\Delta{\mathop{x}\limits_\sim}\!\,_t \,,\,
  E_t[x]+\Delta\widetilde{x}_t ] \;:=\;
  \mathop{\arg\min}_{[a,b]:p_t([a,b])\geq\alpha}(b-a),
\eeqn
where $\Delta\widetilde{x}_t:=\widetilde{x}_t-E_t[x]$
($\Delta{\mathop{x}\limits_\sim}\!\,_t :=
E_t[x]-{\mathop{x}\limits_\sim}\!\,_t$) is the distance from the
right boundary $\widetilde{x}_t$ (left boundary
${\mathop{x}\limits_\sim}\!\,_t$) of the shortest
$\alpha$-credible interval
${\mathop{\widetilde{x}}\limits_\sim}\!\,_t$ to the mean $E_t[x]$
of distribution $p_t$. We can use
$\mathop{\up{\widetilde{x}}}\limits_\simeq \equiv
[\mathop{x}\limits_\simeq , \up{\widetilde{x}}] :=\bigcup_t
{\mathop{\widetilde{x}}\limits_\sim}\!\,_t$ as a (conservative,
but not shortest) robust credible interval, since
$p_t(\mathop{\up{\widetilde{x}}}\limits_\simeq)\geq p_t
({\mathop{\widetilde{x}}\limits_\sim}\!\,_t)\geq\alpha$ for all
$t$. We can upper bound $\up{\widetilde{x}}$ (and similarly lower
bound $\mathop{x}\limits_\simeq$) by
\beq\label{eqCRCI}
  \up{\widetilde{x}} \;=\;
  \max_t(E_t[x]+\Delta\widetilde{x}_t) \;\leq\;
  \max_t E_t[x]+\max_t \Delta\widetilde{x}_t \;=\;
  \up{E[x]}+ \up{\Delta\widetilde{x}}.
\eeq
We have already intensively discussed how to compute upper and
lower quantities, particularly for the upper mean $\up{E[x]}$ for
$x\in\{{\cal F}, {\cal H}, {\cal I}, ... \}$, but the linearization
technique introduced in Section \ref{secAEI} is general enough to
deal with all in $t$ differentiable quantities, including
$\Delta\widetilde{x}_t$. For example for Gaussian $p_t$ with
variances $\sigma_t$ we have
$\Delta\widetilde{x}_t=\kappa\sigma_t$ with $\kappa$ given by
$\alpha=\mbox{erf}(\kappa/\sqrt{2})$, where erf is the error
function (e.g.\ $\kappa=2$ for $\alpha\doteq 95\%)$. We only need to
estimate $\max_t\sigma_t$.

For non-Gaussian distributions, exact expression for
$\Delta\widetilde{x}_t$ are often hard or impossible to obtain and
to deal with. Non-Gaussian distributions depending on some sample
size $n$ are usually close to Gaussian for large $n$ due to the
central limit theorem.
One may simply use $\kappa\sigma_t$ in place of
$\Delta\widetilde{x}_t$ also in this case, keeping in mind that
this could be a non-conservative approximation. More
systematically, simple (and for large $n$ good) upper bounds on
$\Delta\widetilde{x}_t$ can often be obtained and should preferably be used.

Further, we have seen that the variation of sample depending
differentiable functions (like $E_t[x]=E_t[x|\v n]$) w.r.t.\ $t\in\Delta$
are of order ${s\over n+s}$. Since in such cases the standard
deviation $\sigma_t\sim n^{-1/2} \sim \Delta\widetilde{x}_t$ is
itself suppressed, the variation of $\Delta\widetilde{x}_t$ with
$t$ is of order $n^{-3/2}$. If we regard this as negligibly
small, we may simply fix some $t^*\in\Delta$:
\beqn
  \max_t\Delta\widetilde{x}_t = \kappa\sigma_{t^*}+O(n^{-3/2}).
\eeqn
Since $\Delta\widetilde{x}_t$ is ``nearly'' constant, this also
shows that we lose at most $O(n^{-3/2})$ precision in the bound
\req{eqCRCI} (equality holds for $\Delta\widetilde{x}_t$
independent of $t$).

\paradot{Robust credible intervals for mutual information}
Consider the mutual information defined in \req{mi}. The robust
credible interval for $\cal I$ can be estimated as follows.
\beqn
  \up{\widetilde{\cal I}} \;\leq\;
  \up{I}+\up{\Delta\widetilde{\cal I}} \;\leq\;
  I_0+I_R^{ub}+\up{\Delta\widetilde{\cal I}} \;=\;
  I_0+I_R^{ub}+\kappa\sqrt{\Var_{t^*}[{\cal I}]}+O(n^{-3/2}).
\eeqn
Expressions for the variance of $\cal I$ have been
derived in \cite{Hutter:01xentropy}:
\beqn
  \Var_t[{\cal I}] \;=\;
  {1\over n+s}
  \sum_{ij}\u_{ij}\left(\log{\u_{ij}\over
    \u_{i\p}\u_{\p j}}\right)^2 -
  {1\over n+s}\left(\sum_{ij}\u_{ij}\log{\u_{ij}\over
    \u_{i\p}\u_{\p j}}\right)^2 \;+\; O(n^{-2}).
\eeqn
Higher order corrections to the variance and higher moments have
also been derived, but are irrelevant in light of our other
approximations.

\section{Conclusions}\label{secConc}

This is the first work, providing a systematic approach for
deriving closed form expressions for interval estimates for the
Imprecise Dirichlet Model (IDM).
We concentrated on exact and conservative {\em robust} interval
([lower,upper]) estimates for concave functions
$F=\sum_i f_i$ on simplices, like the entropy.
For the conservative estimates we used a first-order Taylor series
expansion in one over the sample size $n$ and bounded the exact
remainder, which widened the intervals by $O(n^{-2})$. This
construction may work for other imprecise models too.
Here is a dilemma, of course: For large $n$ the approximations are
good, whereas for small $n$ the bounds are more interesting, so
the approximations will be most useful for intermediate $n$. More
precise expressions for small $n$ would be highly interesting.
We have also indicated how to
propagate robust estimates from simple functions to composite
functions, like the mutual information.
We argued that a reduced IDM on product spaces, like Bayesian
nets, is more natural and should be preferred in order to improve
predictions. Although improvement is formally only $O(n^{-2})$,
the difference may be significant in Bayes nets or for very small $n$.
Finally, the basics of how to combine robust with credible
intervals have been laid out. Under certain conditions
$O(n^{-3/2})$ approximations can be derived, but the presented
approximations are not conservative.
All in all this work has shown that the IDM has not only interesting
theoretical properties, but that explicit
(exact/conservative/approximate) expressions for robust
(credible) intervals for various quantities can be derived. The
computational complexity of the derived bounds on $F=\sum_i f_i$
is very small, typically one or two evaluations of $F$ or related
functions, like its derivative.
First applications of these (or more precisely, very similar)
results, especially the mutual information, to robust
inference of trees look promising \cite{Hutter:05tree}.

\paradot{Acknowledgements}
I want to thank Peter Walley for introducing the IDM to me, %
Marco Zaffalon for encouraging me to investigate this topic, %
and Jean-Marc Bernard for his feedback on earlier drafts of the paper.

\appendix
\section{Properties of the $\psi$ Function}\label{secPsi}

The digamma function $\psi$ is defined as the logarithmic derivative of
the Gamma function. Integral representations for $\psi$ and its
derivatives are
\beqn
  \psi(z)={d\ln\Gamma(z)\over dz}= {\Gamma'(z)\over\Gamma(z)} =
  \int_0^\infty\left[{e^{-t}\over t}-{e^{-zt}\over
  1-e^{-t}}\right] dt,
  \quad
  \psi^{(k)}(z)=(-1)^{k+1}\int_0^\infty{t^ke^{-zt}\over 1-e^{-t}} dt.
\eeqn
The $h$ function \req{hex} and its first derivative are
\bqan
  h(\u_i) &=& (n_i\!+\!st_i)[\psi(n\!+\!s\!+\!1)-\psi(n_i\!+\!st_i\!+\!1)]/(n\!+\!s),\\
  h'(\u_i) &=& \psi(n\!+\!s\!+\!1)-\psi(n_i\!+\!st_i\!+\!1)-(n_i\!+\!st_i)\psi'(n_i\!+\!st_i\!+\!1),
\eqan
For integral $s$ and at argument $\u_i^0={n_i\over n+s}$ and
$\u_i^0={n_i+s\over n+s}$ we need $\psi$ and $\psi'$ only at integer
values for which the following closed representations exist
\beqn
  \psi(n\!+\!1)=-\gamma+\sum_{i=1}^n{1\over i},\quad
  \psi'(n\!+\!1)={\pi^2\over 6}-\sum_{i=1}^n{1\over i^2},\quad
\eeqn
where $\gamma=0.5772156...$ is Euler's constant. Closed expressions
for half-integer values and fast approximations for arbitrary
arguments also exist. The following asymptotic expansion can be used
if one is interested in $O(({s\over n+s})^2)$ approximations only
(and not rigorous bounds):
\beqn
  \psi(z+1)=\log z + {1\over 2z} - {1\over 12z^2} + O({1\over z^4}),
\eeqn
This shows that $h(\u_i)$ converges to $-\u_i\log\u_i$ for
$n\to\infty$ (and $u_i\to$const.), i.e.\ $H(\vu)$ is close to ${\cal
H}(\vu)$ for large $n$. See \cite[Chp.6]{Abramowitz:74} for details
on the $\psi$ function and its derivatives. From the above
expressions one may show $h''<0$.

\section{Symbols}\label{secSymb}

\begin{tabbing}
  \hspace{2cm} \= \hspace{11cm} \= \kill
  {\bf Symbol }      \> {\bf Explanation}                                                    \\[0.5ex]
  $\delta_{ij}$      \> Kronecker symbol ($\delta_{ij}=1$ for $i=j$ and $\delta_{ij}=0$ for $i\neq j$) \\[0.5ex]
  $\imath,i$         \> Discrete random variable, index/outcome/observation $\in\{1,...,d\}$ \\[0.5ex]
  $d$                \> Dimension of discrete random variable $\imath$                       \\[0.5ex]
  $\pi_i$            \> (Objective/aleatory) probability/chance of $i$                       \\[0.5ex]
  $\log$             \> natural logarithm to basis $e$                                       \\[0.5ex]
  $x_i, \v x, x_\p$  \> Vector $\v x=(x_1,...,x_d),\quad x_\p=x_1+...+x_d,\quad \v x\in\{\v n,\vt,\vu,\vpi,...\}$ \\[0.5ex]
  $t_i, \vt$         \> Initial bias of $i$, bias vector                                     \\[0.5ex]
  $\Deltapi$         \> = $\{\vpi\,:\,\pi_i\geq 0\,\forall i,\, \sum_i\pi_i=1\}\;\;=$ $\vpi$-simplex ($\vpi\in\Deltapi$)  \\[0.5ex]
  $\Delta_{(e)}$     \> = $\{\,\vt\;:\,t_i\geq 0\;\forall i,\; \sum_i t_i\,\smash{\stackrel{(<)}=}\,1\}\;\,\,=$ (extended) $\vt$-simplex ($\vt\in\Delta_{(e)}$)  \\[0.5ex]
  $\Delta'_{(e)}$    \> = $\{\vu\,:\,\u_i\geq\u_i^0\,\forall i,\, \sum_i\u_i\,\smash{\stackrel{(<)}=}\,1\}=$ (extended) $\vu$-simplex ($\vu\in\Delta'_{(e)}$) \\[0.5ex]
  $s$                \> Magnitude of imprecision ($n'_i=st_i$ is virtual observation \#)     \\[0.5ex]
  $\v D$             \> Data/sample $\{i_1,...,i_n\}$                                        \\[0.5ex]
  $n_i,\v n,n$       \> \# of outcomes/observations $i$, \# sample vector, total sample size \\[0.5ex]
  $\delta(\cdot)$    \> Dirac delta distribution $\int f(x)\delta(x)dx=f(0)$                 \\[0.5ex]
  $p(\vpi|\v n)$     \> $\propto\smash{\prod_i\pi_i^{n+st_i-1}}\propto$ Dirichlet posterior              \\[0.5ex]
                     \> (second order/belief/subjective/epistemic probability)               \\[0.5ex]
  $E_\vt[{\cal F}]$  \> Expected value of $\cal F$ w.r.t.\ posterior $p(\vpi|\v n)$          \\[0.5ex]
  w.r.t.             \> with respect to                                                      \\[0.5ex]
  i.i.d.             \> independent and identically distributed                              \\[0.5ex]
  $\u_i^0$           \> $={n_i\over n+s}$                                                    \\[0.5ex]
  $\u_i$             \> $={n_i+st_i\over n+s}=E_\vt[\vpi]$                                   \\[0.5ex]
  $\u_i^*,t_i^*$     \> Origin for Taylor expansion                                          \\[0.5ex]
  $\sigma$           \> $={s\over n+s}=1-\u_\p^0=$ Taylor expansion parameter                \\[0.5ex]
  $O(\sigma^k)$      \> $f(\v n,\vt,s)=O(\sigma^k)$ $\;:\Leftrightarrow\;$ $\exists
                        c\,\forall\v n\in\SetN_0^d,\,\vt\in\Delta,\,s>0$ : $|f(\v n,\v
                        t,s)|\leq c\sigma^k$                                                 \\[0.5ex]
  ${\cal H}(\vpi)$   \> $=-\sum_i\pi_i\log\pi_i=$ entropy of $\vpi$                          \\[0.5ex]
  $H(\vu)$           \> $=\sum_i h(\u_i)=$ expected entropy (see Eq.\req{hex})               \\[0.5ex]
  ${\cal F}(\vpi)$   \> = function of $\vpi$ (${\cal F}\in\{ {\cal H}, {\cal I}, ...\}$)     \\[0.5ex]
  $F(\vu)$           \> = statistic $E_\vt[{\cal F}]$ or general function ($F\in\{H,I,...\}$)\\[0.5ex]
  $F\leqsq G$        \> $:\Leftrightarrow F\leq G \mbox{ and } F=G+O(\sigma^2)$,
                          i.e.\ $G$ is ``good'' upper bound on $F$                           \\[0.5ex]
  $\vu^{\up F},\vt^{\up F}$ \> maximize (and $\vu^{\low F},\vt^{\low F}$ minimize) $F(\vu)$, $\vt\in\Delta$, $\vu\in\Delta'$           \\[0.5ex]
  $\up F$            \> $=\maxo_{\vt\in\Delta} F(\vu)= F(\vu^{\up F})=$ upper value of $F(\vu)$, similarly $\low F$ \\[0.5ex]
  $\up{\low F}$      \> $=[\low F,\up F]=$ robust/Imprecise interval (estimate) of $F$       \\[0.5ex]
  $F_0+F_R(\vu)$     \> $=F(\vu)$ with $F_0=F(\vu^0)$ and $F_R(\vu)=O(\sigma)$               \\[0.5ex]
  $[F_R^{lb},F_R^{ub}]$\> $\supseteq[\low F_R,\up F_R]\ni F_R$ (conservative [lower,upper] bound on $F_R$) \\[0.5ex]
  ${\mathop{\widetilde{F}}\limits_\sim}$ \> $=[{\mathop{F}\limits_\sim}, \widetilde{F}]=$
                         credible interval (estimate) of $F$                                 \\[-0.3ex]
  $\u_{ij},\u_{i\p},\u_{\p j}$\> joint, row, column marginal \\[0.5ex]
  ${\cal I}(\vpi)$   \> $=\sum_{ij}\pi_{ij}\log{\pi_{ij}\over\pi_{i\p}\pi_{\p j}}=$ mutual information of $\vpi$ \\[0.5ex]
  $I(\vu)$           \> $=H(\u_{i\p})+H(\u_{\p j})-H(\u_{ij}) = H_{row}+H_{col}-H_{joint}$   \\[0.5ex]
  $joint,row,col$    \> \hspace{3ex} Index for quantities based on joint, row, column marginal distr. \\[0.5ex]
\end{tabbing}


\begin{small}

\end{small}

\end{document}